\documentclass[amssymb, amstex, amsmath, 11pt]{article}
\newtheorem{theorem}{Theorem}[section]

\newtheorem{corollary}[theorem]{Corollary}

\newtheorem{lemma}[theorem]{Lemma}

\newcommand{\dl}{\par\vspace{5mm}}

\usepackage{mathrsfs}
\usepackage{amssymb}
\begin{document}

\title{ Gr\"obner-Shirshov Basis for the Chinese Monoid\footnote{Supported by the NNSF of China
(No.10771077) and the NSF of Guangdong Province (No.06025062).} }

\author{
{\small Yuqun Chen}\\
 {\small \ School of Mathematical Sciences}\\
{\small \ South China Normal University}\\
{\small \ Guangzhou 510631, China}\\
{\small \ yqchen@scnu.edu.cn}\vspace{4mm}\\
{\small Jianjun Qiu}\footnote {Corresponding author.}\\
{\small \ Mathematics and Computation  Science School }\\
{\small \ Zhanjiang  Normal University}\\
{\small \ Zhanjiang  524048, China}\\
{\small \ jianjunqiu@126.com}
\\
 }

\date{}
 \maketitle

\begin{abstract}\noindent
In this paper, a Gr\"obner-Shirshov basis for  the Chinese monoid is
obtained and an algorithm for the normal form of  the Chinese monoid
is given.
\vspace{4mm}\\
AMS Mathematics Subject Classification (2000): 16S15, 16S35, 20M25.
 \vspace{4mm}\\
Keywords:  Gr\"obner-Shirshov bases; Chinese algebra;  Chinese
monoid; Normal form.
\end{abstract}

\section{Introduction and Preliminaries}
According to the paper \cite{Jc01} (see also \cite{Du94} ), the
Chinese monoid $CH(A,<)$ over a well ordered set $(A, <)$ has the
following presentation
$$
CH(A,<)=sgp\langle A| cba=bca=cab,\ a,b,c \in A,a\leq b\leq c
\rangle.
$$
Let $K$ be a field. The semigroup algebra $K\langle A|T\rangle$ $ $
of  the  Chinese monoid  $CH(A)$ is called the Chinese algebra.

In the cited paper  \cite{Jc01}, some special tools were invented to
study the word problem and other enumeration problems for the
Chinese monoid. These are first of all the Chinese staircase and the
insertion algorithm. Using these tools, the authors in \cite{Jc01}
were able to construct normal form of words of the Chinese monoid
and the algorithm to find the normal form of a word.

In the present paper, we are going to simplify some part of the
paper  \cite{Jc01} by  using the Gr\"obner-Shirshov bases theory for
associative algebras. The main new result of this  paper is the
Gr\"obner-Shirshov basis for the Chinese monoid. As a simple
corollary, we get  a normal form of words that coincide with the
staircase normal form found  in the paper \cite{Jc01}. Also the
algorithm of elimination of leading words of relations gives rise
the insertion algorithm of the same paper.

In the rest of  this section, we recall the definitions of  the
Gr\"obner-Shirshov bases and Composition-Diamond Lemma for
associative algebra.

Let $K$ be a field, $K\langle X\rangle$ the free associative algebra
over $K$ generated by $X$. Denote $X^*$ the free monoid generated by
$X$, where the empty word is the identity which is denoted by 1. For
a word $w\in X^*$, we denote the length of $w$ by $|w|$. Let $X^*$
be a well ordered set. Then every nonzero  polynomial  $f\in
K\langle X\rangle$ has the leading word $\bar{f}$. If the
coefficient of $\bar{f}$ in $f$ is equal to 1, then  $f$ is callled
monic.

Let $f$ and $g$ be two monic polynomials in $K\langle X\rangle$.
Then, there are two kinds of compositions:

$(1)$ If \ $w$ is a word such that $w=\bar{f}b=a\bar{g}$ for some
$a,b\in X^*$ with $|\bar{f}|+|\bar{g}|>|w|$, then the polynomial
 $(f,g)_w=fb-ag$ is called the intersection composition of $f$ and
$g$ with respect to $w$. The word $w$ is called an ambiguity of
intersection.

$(2)$ If  $w=\bar{f}=a\bar{g}b$ for some $a,b\in X^*$, then the
polynomial $(f,g)_w=f - agb$ is called the inclusion composition of
$f$ and $g$ with respect to $w$. The word $w$ is called an ambiguity
of inclusion.

If $g$ is monic, $\bar{f}=a\bar{g}b$  and $\alpha$ is the
coefficient of the leading term $\bar{f}$,  then transformation
$f\mapsto f-\alpha agb$ is called  elimination of the leading word
(ELW) of $g$ in $f$.

Let $S\subseteq$ $K\langle X\rangle$ with each $s\in S$ monic. Then
the composition $(f,g)_w$ is called trivial modulo $(S,\ w)$ if
$(f,g)_w=\sum\alpha_i a_i s_i b_i$, where each $\alpha_i\in k$,
$a_i,b_i\in X^{*}, \ s_i\in S$ and $a_i \overline{ s_i }b_i<w$. If
this is the case, then we write
$$
(f,g)_w\equiv0\quad mod(S,w).
$$
In general, for $p,q\in K\langle X\rangle$, we write
$$
p\equiv q\quad mod (S,w)
$$
which means that $p-q=\sum\alpha_i a_i s_i b_i $, where each
$\alpha_i\in K,a_i,b_i\in X^{*}, \ s_i\in S$ and $a_i \overline{
s_i} b_i<w$.

We call the set $S$ endowed with the well order $<$ a
Gr\"{o}bner-Shirshov basis  for  $K\langle X|S\rangle$ if any
composition of polynomials in $S$ is trivial modulo $S$.

A well order $<$ on $X^*$ is monomial if for $u, v\in X^*$, we have
$$
u < v \Rightarrow w_{1}uw_{2} < w_{1}vw_{2},  \ for \  all \
 w_{1}, \ w_{2}\in  X^*.
$$

The following lemma was proved by Shirshov \cite{Sh} for  free Lie
algebras (with deg-lex ordering) in 1962 (see also Bokut
\cite{b72}). In 1976, Bokut \cite{b76} specialized the approach of
Shirshov to associative algebras (see also Bergman \cite{b}). For
commutative polynomials, this lemma is known as the Buchberger's
Theorem (see \cite{bu65} and \cite{bu70}).

\ \

\noindent{\bf Composition-Diamond Lemma.} \ Let $K$ be a field,
$\mathscr{A}=K \langle X|S\rangle=K\langle X\rangle/Id(S)$ and $>$ a
monomial order on $X^*$, where $Id(S)$ is the ideal of $K \langle
X\rangle$ generated by $S$. Then the following statements are
equivalent:
\begin{enumerate}
\item[(1)] $S $ is a Gr\"{o}bner-Shirshov basis.
\item[(2)] $f\in Id(S)\Rightarrow \bar{f}=a\bar{s}b$
for some $s\in S$ and $a,b\in  X^*$.
\item[(3)] $Irr(S) = \{ u \in X^* |  u \neq a\bar{s}b ,s\in S,a ,b \in X^*\}$
is a basis of the algebra $\mathscr{A}=K\langle X | S \rangle$.
\end{enumerate}

If a subset $S$ of $K\langle X\rangle$ is not a Gr\"{o}bner-Shirshov
basis, then we can add to $S$ all nontrivial compositions of
polynomials of $S$, and by continuing this process (maybe
infinitely) many times, we eventually obtain a Gr\"{o}bner-Shirshov
basis $S^{comp}$. Such a process is called the Shirshov algorithm.

If $S$ is a set of ``semigroup relations" (that is, the polynomials
of the form $u-v$, where $u,v\in X^{*}$), then any nontrivial
composition will have the same form. As a result, the set $S^{comp}$
also consists of semigroup relations.

Let $M=sgp\langle X|S\rangle$ be a semigroup presentation. Then $S$
is a subset of $K\langle X\rangle$  and hence one can find a
Gr\"{o}bner-Shirshov basis $S^{comp}$. The last set does not depend
on $K$, and as mentioned before, it consists of semigroup relations.
We will call $S^{comp}$ a Gr\"{o}bner-Shirshov basis of $M$. This is
the same as a Gr\"{o}bner-Shirshov basis of the semigroup algebra
$KM=K\langle X|S\rangle$.  If $S$ is  a Gr\"{o}bner-Shirshov basis
of the semigroup  $M=sgp\langle X|S\rangle$, then $Irr(S)$ is a
normal form for $M$.

\section{Chinese Algebras and Monoids }
In this section, we find a Gr\"{o}bner-Shirshov basis of the Chinese
monoid  (Chinese algebra) and give an algorithm to find the normal
form of words of the Chinese monoid.

Let $A=\{x_i|i\in I\}$ be a well ordered set. The Chinese congruence
is the congruence  on $A^*$ generated by $T$, where $T$ consists of
the relations:
\begin{eqnarray*}
&& x_ix_jx_k=x_ix_kx_j= x_jx_ix_k\ \ for\  every \ i>j>k,\\
 &&x_ix_jx_j= x_jx_ix_j,\ x_ix_ix_j= x_ix_jx_i \ \ for\
every \ i>j.
\end{eqnarray*}
The Chinese monoid $CH(A)$ is the quotient monoid of the free monoid
$A^*$ by the Chinese congruence, i.e.,  $CH(A) =sgp\langle
A|T\rangle$.  Let $K$ be a field. The semigroup algebra $K\langle
A|T\rangle$ $ $ of  the  Chinese monoid  $CH(A)$ is called the
Chinese algebra.

Now, we define an order on $A^*$ by the deg-lex order.

Let $S$ be the set which  consists of  the following polynomials:
\begin{enumerate}
\item[1] \ $x_ix_jx_k-x_jx_ix_k$,
\item[2] \
$x_ix_kx_j-x_jx_ix_k$,
 \item[3] \
 $x_ix_jx_j-x_jx_ix_j$,
 \item[4] \
$x_ix_ix_j-x_ix_jx_i $,
\item[5] \ $x_ix_jx_ix_k-x_ix_kx_ix_j$,
\end{enumerate}
where $x_i, x_j, x_k \in A,\ i>j>k$.

 It is easy to check  that  $ sgp\langle A|T\rangle=sgp\langle
A|S\rangle$.  In the following theorem, we shall prove that $S$ is a
Gr\"{o}bner-Shirshov basis of the Chinese monoid $CH(A)$.

\begin{theorem} \label{2.3}
With the deg-lex order on $A^*$, $S$ is a Gr\"{o}bner-Shirshov basis
of the Chinese monoid $CH(A)$.
\end{theorem}
\noindent {\bf Proof.}  Denote by $i\wedge j$ the composition of the
polynomials of type $i$ and type $j$ in $S$. For example,
\begin{eqnarray*}
1\wedge4&=&(x_ix_jx_k-x_jx_ix_k,x_kx_kx_{j_{1}}-x_kx_{j_{1}}x_k
)_{x_ix_jx_kx_kx_{j_1}}\\
&=&
(x_ix_jx_k-x_jx_ix_k)x_kx_{j_{1}}-x_ix_j(x_kx_kx_{j_{1}}-x_kx_{j_{1}}x_k),
\end{eqnarray*}
where $ i>j>k>j_1$.

The possible ambiguities of  two polynomials in $S$ are only as
below (instead of $x_i$, we write just $i$). In the following list,
$i\wedge j\ \ \ w$ means that $w$ is the  ambiguity   of  the
composition $i\wedge j$.
\begin{tabbing}
 $1\wedge1\ \ \ ijkk_1,\ i>j>k>k_1$  \qquad  \qquad \qquad \= $1\wedge1\ \ \ ijkj_1k_1, \
 i>j>k>j_1>k_1$\\ [0.5ex]
 $1\wedge 2\ \ \ ijj_1k, \  i>j>j_1>k$ \>$ 1\wedge 2\ \ \ ijkj_1k_1, \
 i>j>k>j_1>k_1$\\ [0.5ex]
$1\wedge3 \ \ \ ijkk, \ i>j>k$  \> $1\wedge3\ \ \ ijk{j_1}{j_1}, \
i>j>k>j_1$\\ [0.5ex]
$1\wedge4\ \ \  ijkk{j_1}, \ i>j>k>j_1$  \> $ 1\wedge5  \ \ \  ijkj{k_1},\  i>j>k>k_1$ \\[0.5ex]
$ 1\wedge5\ \ \ ijk{j_1}k{k_1}, \ i>j>k>j_1>k_1$ \> $ 2\wedge 1\ \ \
ikj{j_1}{k_1},\ i>j>k,\ j>j_1>k_1$\\ [0.5ex]
$2\wedge 2\ \ \ ikj{k_1}{j_1},\ i>j>k,\ j>j_1>k_1$ \> $2\wedge 3\ \ \ ikj{j_1}{j_1},\ i>j>k,\  j>j_1$\\[0.5ex]
$2\wedge 4\ \ \ ikjj{j_1},\ i>j>k,\  j>j_1$\> $2\wedge5\ \ \ ikj{j_1}j{k_1}, \ i>j>k,\ j>j_1>k_1$\\[0.5ex]
$3\wedge1 \ \ \ ijj{j_1}{k_1}, \ i>j>j_1>k_1$ \>$3\wedge2\ \ \ ijj{k_1}{j_1}, \  i>j>j_1>k_1$\\[0.5ex]
$3\wedge3\ \ \ ijj{j_1}{j_1},\ i>j>j_1$ \> $3\wedge4\ \ \ ijj{j_1}, \ i>j>j_1$\\[0.5ex]
$  3\wedge5\ \ \ ijj{j_1}j{k_1},\  i>j>j_1>k_1$ \> $ 4\wedge1\ \ \ iij{k}, \ i>j>k$\\[0.5ex]
$ 4\wedge1\ \ \ iij{j_1}{k_1}, \ i>j>j_1>k_1$\>$4\wedge2\ \ \ iij{j_1},\ i>j_1>j$\\[0.5ex]
$4\wedge2\ \ \ iij{k_1}{j_1}, \ i>j>j_1>k_1$\> $ 4\wedge3\ \ \ iijj,\ i>j$\\[0.5ex]
$ 4\wedge3\ \ \ iij{j_1}{j_1},\  i>j>j_1$\>$4\wedge4\ \ \ iijj{j_1},\ i>j>j_1$\\[0.5ex]
$ 4\wedge5 \ \ \ iiji{j_1},\ i>j>j_1$\>$4\wedge5\ \ \ iij{k_1}j{j_1},\ i>j>j_1>k_1$\\[0.5ex]
$5\wedge1\ \ \ ijik{k_1},\  i>j>k>k_1$\>$5\wedge1\ \ \ ijik{j_1}{k_1},\ i>j>k>j_1>k_1$\\[0.5ex]
$5\wedge2\ \ \ ijik{j_1},\ i>j>k,\ i>j_1>k$ \>$ 5\wedge2\ \ \ ijik{k_1}{j_1},\ \ i>j>k>j_1>k_1$\\[0.5ex]
$5\wedge3\ \ \ ijikk, \ \ i>j>k$\>$ 5\wedge3\ \ \ ijik{j_1}{j_1},\  i>j>k>j_1$\\[0.5ex]
$ 5\wedge4\ \ \ ijikk{j_1},\ i>j>k>j_1$ \>$5\wedge5\ \ \ ijiki{k_1},\  i>j>k>k_1$\\[0.5ex]
$  5\wedge5\ \ \ ijik{j_1}k{k_1},\  i>j>k>k_1$
\end{tabbing}
It is easy to  check that all these compositions  are trivial modulo
$(S,w)$. Here, for example, we just check $3\wedge2$.
\begin{eqnarray*}
3\wedge2&=&(x_ix_jx_j-x_jx_ix_j,x_jx_{k_1}x_{j_1}-x_{j_1}x_jx_{k_1})_{x_ix_jx_jx_{k_1}x_{j_1}}\\
&=&-x_jx_ix_jx_{k_1}x_{j_1}+x_ix_jx_{j_1}x_jx_{k_1}\\
&\equiv&-x_jx_jx_ix_{k_1}x_{j_1}+x_jx_ix_{j_1}x_jx_{k_1}\\
&\equiv&-x_jx_jx_{j_1}x_ix_{k_1}+x_jx_jx_ix_{j_1}x_{k_1}\\
&\equiv&-x_jx_jx_{j_1}x_ix_{k_1}+x_jx_jx_{j_1}x_ix_{k_1}\\
&=&0,
\end{eqnarray*}
where $i>j>j_1>k_1$.

So  $S$ is a Gr\"{o}bner-Shirshov basis of the Chinese monoid $CH(A)$. \hfill $\blacksquare$\\

Let $\Omega$ be the set which consists of words on $A$ of the form
$u_n=w_1w_2\cdots w_n, n\geq0$, where $ w_1 = x_1^{t_{1,1}}, \
w_k=(\Pi_{1\leq i \leq k-1} (x_kx_i)^{t_{k,i}})x_{k}^{t_{k,k}}, \
2\leq k\leq n, \ x_i\in A, \ x_1<x_2< \cdots <x_n$ and all exponents
are non-negative integers.

By using the ELW, we have the following algorithm (the insertion
algorithm
in \cite{Jc01}).\\

 \noindent  \textbf{Algorithm} \ \  Let $u_n=w_1w_2\cdots
w_n \in \Omega, $   $x\in A $  and  $N(w_n)=\{i|t_{n, i}\neq 0\}$.
Define $u_n\star x$ by induction on $n$:
\begin{enumerate}
\item[(1)] \  If $x>x_n$, then $u_n\star x= u_nw_{n+1}$, where
$w_{n+1}=x$.

\item[(2)] \  If $x=x_n$, then $u_n\star  x= u_{n-1} w_n^{'}$,
where $w_n^{'}=w_n|_{t_{n,n}\rightarrow t_{{n,n}}+1}$.

\item[(3)] \  If $x_1\leq x=x_k<x_n$, let  $i$ be the greatest
number of $N(w_n)$, where $i=k$ if $N(w_n)=\varnothing$. There are
three cases to consider:
\begin{enumerate}
\item[(a)] \  If $x\geq x_i$ , then $u_n\star  x= (u_{n-1}\star
    x)w_n$.
\item[(b)] \  If $x< x_i< x_n$, then  $u_n\star  x= (u_{n-1}\star
    x_i)w_n^{'}$, where $w_n^{'}=w_n|_{t_{n,i}\rightarrow t_{n, i}-1,\  t_{n, k}\rightarrow
    t_{n, k}+1}$.
\item[(c)] \  If $x< x_i= x_n$, then $u_n\star  x= u_{n-1}
  w_n^{'}$, where $w_n^{'}=w_n|_{t_{n, n}\rightarrow t_{n, n}-1,\  t_{n, k}\rightarrow
    t_{n, k}+1}$.
\end{enumerate}
\item[(4)] \  If $x<x_1$,
 then $u_n \star  x= u_{n+1} \star  x$, where $x_0=x$, $w_0=x_0^0$ and
$u_{n+1}=w_0u_n$.
\end{enumerate}

\begin{lemma} \label{2.4}
 Let $u_n=w_1w_2\cdots
w_n \in \Omega$  and  $x\in A $. Then  $u_n\star x\in \Omega
 $ and in $ K\langle
A|S\rangle, \ u_nx = u_n\star  x$.
\end{lemma}
\noindent {\bf Proof.} We prove the results by induction on $n$.

If  $n=0$ or 1, then by using
  the ELW of the relations in $S$, it is easy to check that
$u_1\star  x\in \Omega
 $ and $u_1x = u_1\star  x$ in $ K\langle
A|S\rangle$.

 Now we assume $n\geq 2$.  For the cases of (1), (2) and (3) (c), by
  the ELW, we have the results.

For (3) (a),  we have  $u_{n-1}\star  x\in \Omega
 $ and $u_{n-1}x = u_{n-1}\star  x$ in $K\langle
A|S\rangle$ by induction. Since $x\geq x_i$,  $xw_n=w_nx$ in
$K\langle A|S\rangle$ by the ELW. It follows that $u_n\star  x\in
\Omega $ and $u_nx=u_{n-1}w_nx=u_{n-1}xw_n=(u_{n-1}\star
x)w_n=u_n\star  x$ in $ K\langle A|S\rangle$.

For (3) (b), by using the ELW,  we have $x_iw_n^{'}=w_nx$ in
$K\langle A|S\rangle$ for $x< x_i< x_n$. So, by induction,
$u_{n}\star x\in \Omega $ and
$u_nx=u_{n-1}w_nx=u_{n-1}x_iw_n^{'}=(u_{n-1}\star
x_i)w_n^{'}=u_n\star  x$ in $ K\langle A|S\rangle$.

(4) follows from (3). \hfill $\blacksquare$

 \begin{theorem} \label{2.5}
The set  $\Omega$ is a basis of the Chinese algebra $K\langle
A|S\rangle$ as a $K$-space.
 \end{theorem}
 \noindent {\bf Proof.}
 By Theorem \ref{2.3} and the  Composition-Diamond Lemma,
 $ Irr(S)=\{ w\in A^*|w\neq a\bar{s}b ,s\in S,a ,b \in A^*\}$ is a
 linearly
 basis of $K\langle
A|S\rangle$. Clearly, $\Omega \subseteq Irr(S)$ whence $\Omega$ is
linearly  independent. For any word $w\in Irr(S)$, we use induction
on $|w|$ to prove that there exists a $t\in \Omega $ such that $w=t$
in $ K\langle A|S\rangle$.
  If $|w|=1$, then $w\in \Omega $.
Assume that  $|w|\geq 2$ and let $w=w_1x$. Then, by  induction and
Lemma \ref{2.4}, we have $t_1\in
  \Omega$  such that $w_1=t_1$ in $ K\langle A|S\rangle$  and $t_1\star  x\in \Omega$ whence  $w=w_1x= t_1x=t_1\star
  x$ in $ K\langle
A|S\rangle$. So $\Omega$ is a linearly basis of $ K\langle
A|S\rangle$, i.e., $\Omega = Irr(S)$.\hfill $\blacksquare$

\begin{corollary} {\em (\cite{Jc01}  Theorem 2.1)}
The  Chinese monoid  $CH(A)$ has a normal  form of words
$u_n=w_1w_2\cdots w_n,\ n\geq 0$ {\em (}the row normal form of a
staircase{\em )}.
\end{corollary}

 \noindent {\bf Proof.} By  Theorem \ref{2.3} and  Theorem \ref{2.5}.\hfill
$\blacksquare$

\ \

\noindent{\bf Acknowledgement}: The authors would like to express
their deepest gratitude to Professor L. A. Bokut for his kind
guidance, useful discussions and enthusiastic encouragement during
his visit to the South China Normal University.

\end{document}